\documentclass[reqno,12pt]{amsart}
\usepackage{a4wide}
\usepackage[ps,all]{xy}
\usepackage[dvipsnames]{xcolor}
\usepackage[pdfusetitle]{hyperref}
\hypersetup{
 colorlinks=true,
 linkcolor=Blue,
 citecolor=Blue,
 urlcolor=Blue,
 }
 \usepackage[alphabetic,msc-links,initials]{amsrefs}
\usepackage{graphicx}
\usepackage{mathptmx}
\usepackage{latexsym}
\usepackage{amsmath}
\usepackage{mathtools,amssymb}
\usepackage{tikz}
\usepackage{enumerate}
\usepackage[ignoredisplayed]{enumitem}
\usetikzlibrary{arrows}
\usepackage{newtxtext}
\usepackage[varg]{newtxmath}
\theoremstyle{definition}
\newtheorem{theorem}{Theorem}
\newtheorem*{theorem*}{Theorem}
\newtheorem{definition}[theorem]{Definition}
\newtheorem*{definition*}{Definition}
\newtheorem{proposition}[theorem]{Proposition}
\newtheorem*{proposition*}{Proposition}
\newtheorem{lemma}[theorem]{Lemma}
\newtheorem*{lemma*}{Lemma}
\newtheorem{def-thm}[theorem]{Definition-Theorem}
\newtheorem*{def-thm*}{Definition-Theorem}
\newtheorem{def-pro}[theorem]{Definition-Proposition}
\newtheorem*{def-pro*}{Definition-Proposition}
\newtheorem{corollary}[theorem]{Corollary}
\newtheorem*{corollary*}{Corollary}
\newtheorem{example}[theorem]{Example}
\newtheorem*{example*}{Example}
\newtheorem{remark}[theorem]{Remark}
\newtheorem*{remark*}{Remark}

\numberwithin{theorem}{section}
\numberwithin{equation}{section}
\newcommand{\sgn}[2]{\operatorname{sgn}(#1, #2)}
\newcommand{\eps}[1]{{\varepsilon}_{#1}}
\newcommand{\gcdd}[2]{\operatorname{gcd}(#1, #2)}
\newcommand{\diag}[1]{\operatorname{\mathrm{diag}}(#1)}

\newcommand{\length}[1]{\mathit{l}(#1)}
\newcommand{\summand}[1]{| #1 |}
\newcommand{\D}[1]{{\mathbb{D}} #1}
\newcommand{\Hom}[3]{\operatorname{\mathrm{Hom}}_{#1}(#2, #3)}

\newcommand{\ann}[1]{\operatorname{ann}#1}
\newcommand{\Ext}[4]{\operatorname{\mathrm{Ext}}^{#1}_{#2}(#3, #4)}
\newcommand{\Tor}[4]{\operatorname{Tor}_{#1}^{#2}(#3, #4)}

\newcommand{\soc}[1]{\operatorname{\mathrm{soc}}#1}
\newcommand{\add}[1]{\operatorname{\mathsf{add}}#1}
\newcommand{\Fac}[1]{\operatorname{\mathsf{Fac}}#1}
\newcommand{\Sub}[1]{\operatorname{\mathsf{Sub}}#1}
\newcommand{\Rep}[1]{\operatorname{\mathsf{Rep}}#1}
\newcommand{\istrep}[1]{\operatorname{\overline{\mathsf{rep}}}#1}
\newcommand{\pstrep}[1]{\operatorname{\underline{\mathsf{rep}}}#1}
\newcommand{\rep}[1]{\operatorname{\mathsf{rep}}#1}
\newcommand{\tors}[1]{\operatorname{\mathsf{tors}}#1}
\newcommand{\torf}[1]{\operatorname{\mathsf{torf}}#1}
\newcommand{\ftors}[1]{\operatorname{\mathsf{f\mathchar`-tors}}#1}
\newcommand{\ftorf}[1]{\operatorname{\mathsf{f\mathchar`-torf}}#1}
\newcommand{\Ker}[1]{\operatorname{\mathrm{Ker}}#1}
\newcommand{\Img}[1]{\operatorname{\mathrm{Im}}#1}

\newcommand{\fgpre}[2]{{\Pi}{\,}(#1, #2)}
\newcommand{\affPi}{\widetilde{\Pi}}
\newcommand{\affI}[1]{\widehat{I}_{#1}}
\newcommand{\tildeI}[1]{\widetilde{I}_{#1}}
\newcommand{\lfRep}[1]{\operatorname{\mathsf{Rep}_{\mathrm{l\mathchar`.f\mathchar`.}}}#1}
\newcommand{\lfrep}[1]{\operatorname{\mathsf{rep}_{\mathrm{l\mathchar`.f\mathchar`.}}}#1}
\newcommand{\proj}[1]{\operatorname{\mathsf{proj}}#1}

\newcommand{\itrigid}[1]{\mathsf{i\tau \mathchar`-rigid}{\,}#1}
\newcommand{\sttilt}[1]{\mathsf{s\tau \mathchar`-tilt}{\,}#1}
\newcommand{\stmtilt}[1]{\mathsf{s\tau^- \mathchar`-tilt}{\,}#1}
\newcommand{\mirr}[1]{\mathsf{m\mathchar`-Irr}{\,}#1}
\newcommand{\jirr}[1]{\mathsf{j\mathchar`-Irr}{\,}#1}

\newcommand{\pdim}[2]{\mathrm{proj.dim}_{#1}{\,}#2}

\newcommand{\op}[1]{{#1}^\mathrm{op}}
\newcommand{\tenso}[1]{{\otimes}_{#1}}

\newenvironment{enumerate2}
{
\begin{list}{\rm{(\arabic{enumi}).}}
{
\usecounter{enumi}
\setlength{\topsep}{0em}
\setlength{\itemindent}{0em}
\setlength{\leftmargin}{2.2em}
\setlength{\rightmargin}{0em}
\setlength{\labelsep}{0.3em}
\setlength{\labelwidth}{3em}
\setlength{\itemsep}{0em} 
\setlength{\parsep}{0em} 
\setlength{\listparindent}{1em} 
}
}{
\end{list}
}

\newenvironment{itemize2}
{
\begin{list}{\rm{(\roman{enumi}).}}
{
\usecounter{enumi}
\setlength{\topsep}{0em}
\setlength{\itemindent}{0em}
\setlength{\leftmargin}{2.2em}
\setlength{\rightmargin}{0em}
\setlength{\labelsep}{0.3em}
\setlength{\labelwidth}{3em}
\setlength{\itemsep}{0em}
\setlength{\parsep}{0em}
\setlength{\listparindent}{1em}
}
}{
\end{list}
}
\title[Module categories of generalized preprojective algebras]{On the module categories of generalized preprojective algebras of Dynkin type}
\author[K. Murakami]{Kota Murakami}
\address[K. Murakami]{Department of Mathematics, Kyoto University, Kitashirakawa Oiwake-cho, Sakyo-ku, 
Kyoto 606-8502, Japan}
\email{k-murakami@math.kyoto-u.ac.jp}
\date{\today}
\begin{document}
\maketitle
\begin{abstract}
For a symmetrizable GCM $C$ and its symmetrizer $D$, Geiss-Leclerc-Schr\"oer [Invent. Math. 209 (2017)] has introduced a generalized preprojective algebra $\Pi$ associated to $C$ and $D$, that contains a class of modules, called locally free modules. We show that any basic support $\tau$-tilting $\Pi$-module is locally free and gives a classification theorem of torsion-free classes in $\rep{\Pi}$ as the generalization of the work of Mizuno [Math. Z. 277 (2014)].
\end{abstract}
\section{Introduction}
\label{intro}
In the progress of representation theory of quivers, Gabriel \cite{MR0332887} has shown that a connected quiver has finitely many indecomposable representations if and only if the underlying graphs are Dynkin diagrams of $A, D, E$ of finite types. In particular, the isoclasses of indecomposable representations corresponds to its positive root systems bijectively via their dimension vectors. Later, Kac \cites{MR557581,MR677715} has extended this description from the cases of symmetric Dynkin types to the case containing any symmetric affine types. That is, it has given a characterization of indecomposable representations of the quivers in terms of positive real roots and positive imaginary roots. On the other hand, Gel'fand-Ponomarev \cite{MR545362} has introduced preprojective algebras for acyclic quivers in order to develop Auslander-Reiten theory for the quivers.  In the works of Buan-Iyama-Reiten-Scott \cite{MR2521253} and Mizuno \cite{Mizuno2014}, tilting theory of the preprojective algebras for an acyclic quiver is studied in terms of the corresponding Weyl group $W$ (Iyama-Reading-Reiten-Thomas \cite{Iyama2016}, Asai \cite{Asai2018}, Kimura \cite{Kimura2018}). In their study, many notions in tilting theory of preprojective algebras are studied in terms of certain collection of two-sided ideals $I_w$ corresponding to the $w\in W$. This also has some applications to algebraic Lie theory (e.g. Baumann-Kamnitzer-Tingley \cite{MR3270589}, Gei\ss-Leclerc-Schr\"oer \cite{MR2822235}).

In search of its generalization to symmetrizable Kac-Moody algebras, Geiss-Leclerc-Schr\"oer \cite{MR3660306} has introduced a class of $1$-Iwanaga-Gorenstein algebras attached to a generalized Cartan matrix (GCM) and its symmetrizer. It has also introduced generalized preprojective algebras associated to the $1$-Iwanaga-Gorenstein algebra and has developed a class of modules, that is, locally free modules. These algebras and its locally free module categories share many features with the classical study of preprojective algebras with symmetric GCMs replaced by symmetrizable ones. Recently, Fu-Geng \cite{Fu2018} has described some properties of the two-sided ideals $I_w$ ($w \in W$) of generalized preprojective algebras. The goal of this paper is to develop the relationship between tilting theory of generalized preprojective algebras and locally free modules following the work of Fu-Geng \cite{Fu2018}.

Let $C$ be a symmetrizable GCM and let $D$ be a symmetrizer of $C$. We denote a generalized preprojective algebra associated with $C$ and $D$ by $\Pi=\fgpre{C}{D}$. Let $W$ be the Weyl group of Kac-Moody Lie algebra associated with $C$.
\begin{theorem}[$\doteq$ Theorem \ref{I_w-locallyfree}]\label{main:1}
	Let $C$ be a symmetrizable GCM of Dynkin type and let $D$ be any symmetrizer of $C$. For each $w\in W$, the two sided ideal $I_w$ of $\fgpre{C}{D}$ is locally free. In particular, any basic support $\tau$-tilting $\Pi$-modules are locally free. 
\end{theorem}
\begin{theorem}[$\doteq$ Theorem \ref{dualityofI}, Corollary \ref{Pitorstorf}]\label{main:2}
	Let $C$ be a symmetrizable GCM of Dynkin type and let $D$ be any symmetrizer of $C$. We have the following two bijections.
	\begin{align*}
    W &\longrightarrow \stmtilt{\Pi}
    &W &\longrightarrow \torf{\Pi}\\
    w &\longmapsto \Pi/I_w,
	&w &\longmapsto \Sub{(\Pi / I_{w})}.
    \end{align*}
    Here, $\stmtilt{\Pi}$ (resp. $\torf{\Pi}$) is the set of isoclasses of basic support $\tau^-$-tilting modules in the sense of Adachi-Iyama-Reiten \cite{MR3187626} (resp. the set of torsion-free classes).
\end{theorem}
This paper is organized as follows. In section \ref{Prelimi}, we first review the definition of algebras $\fgpre{C}{D}$ by Geiss-Leclerc-Schr\"oer \cite{MR3660306} and some basic properties. Then, we review its $\tau$-tilting theory and the work of Fu-Geng \cite{Fu2018}. In section \ref{we}, we prove Theorem \ref{main:1} and \ref{main:2}.
\section{Preliminaries}\label{Prelimi} 
Throughout this paper, $K$ denotes an arbitrary field, and a $K$-algebra always means an associative algebra with a unit over $K$. For an algebra $\Lambda$, its $\Lambda$-module means a left $\Lambda$-module. We denote the module category over $\Lambda$ by $\Rep{\Lambda}$, and the full subcategory of finite dimensional modules by $\rep{\Lambda}$. We denote the full subcategory of finitely generated projective modules by $\proj{\Lambda}$. We refer to Assem-Simson-Skowro\'nski \cite{MR2197389} for basic concepts and terminologies about algebras and quiver representations. For a lattice $(L, \leq)$, we denote the set of join-irreducible (resp. meet-irreducible) elements by $\jirr{L}\coloneqq\{j\in L\mid \text{$j$ covers a unique element in $L$}\}$ (resp. $\mirr{L}\coloneqq\{m\in L\mid \text{$m$ is covered by a unique element in $L$}\}$). The terminologies about lattice theory appeared in this paper can be found in Iyama-Reading-Reiten-Thomas \cite{Iyama2016}.
\subsection{Preprojective algebras associated with symmetrizable Cartan matrices}
In this subsection, we review the representation theory of preprojective algebras with symmetrizable Cartan matrices introduced in  Geiss-Leclerc-Schr\"oer \cite{MR3660306}. 
\begin{definition}\label{GCM}
	A matrix $C=(c_{ij}) \in M_n(\mathbb{Z})$ is called a \textit{generalized Cartan matrix} (GCM), if it satisfies the following three conditions:
	\begin{enumerate2}
		\item[(C1)] $c_{ii} =2$ for any $1\leq i\leq n$;
		\item[(C2)] If $i \neq j$, then $c_{ij} \leq 0$;
		\item[(C3)] $c_{ij} \neq 0$ if and only if $c_{ji} \neq 0$.
	\end{enumerate2}	
	In particular, if $C=(c_{ij})$ satisfies the following (C4), then $C$ is called \textit{symmetrizable}:
	\begin{enumerate2}
		\item[(C4)] There is a diagonal matrix $D=\diag{c_1, \dots , c_n}{\,}(c_i \in \mathbb{Z}, c_i \geq 1)$ such that $DC$ is a symmetric matrix.
	\end{enumerate2}	
In the condition (C4), the matrix $D$ is called a \textit{symmetrizer} of $C$.
\end{definition}
The following quadratic forms $q_{C}$ and graphs $\Gamma(C)$ give a classification of GCMs:
\begin{definition}
Let $C=(c_{ij}) \in M_n(\mathbb{Z})$ be a symmetrizable GCM and $D=\diag{c_1, \dots , c_n}$ be a symmetrizer of $C$.
\begin{enumerate2}
	\item The graph $\Gamma(C)$ has vertices $1, \dots , n$. An edge between $i$ and $j$ exists in $\Gamma(C)$ if and only if $c_{ij}<0$. The edge has its value $(|c_{ji}|, |c_{ij}|)$: 
$$\begin{tikzpicture}[auto]
	\node (a) at (0, 0) {$i$}; \node (b) at (2.4, 0) {$j$}; 
	\draw[-] (a) to node {$\scriptstyle (|c_{ji}|, |c_{ij}|)$} (b);
\end{tikzpicture}.$$
We call $\Gamma(C)$ the \textit{valued graph} of $C$. If $\Gamma(C)$ is a connected graph, then $C$ is called \textit{connected}.
\item We define the form $q_C \colon\mathbb{Z}^n\longrightarrow \mathbb{Z}$ by
      $$q_C =\sum_{i=1}^n c_i X_i^2 -\sum_{i<j}c_i |c_{ij}|X_iX_j.$$
     Since we have $c_i|c_{ij}| =c_j|c_{ji}|$, the form $q_{C}$ is symmetric. If $q_C$ is positive definite (resp. positive semi-definite), then $C$ is called \textit{Dynkin type} (resp. \textit{Euclidean type}).
\end{enumerate2}
\end{definition}
\begin{remark}
	If $C$ is a connected symmetrizable GCM, there exists a unique \textit{minimal} symmetrizer. That is, any symmetrizer $D'$ of $C$ is equal to $mD$ for the minimal symmetrizer $D$ and some positive integer $m$.
\end{remark}
In Geiss-Leclerc-Schr\"oer \cite{MR3660306}, a generalization of the representation theory of acyclic quivers and their preprojective algebras is given by the data of symmetrizable GCMs and their symmetrizers. First, we define a quiver from a symmetrizable GCM.
\begin{definition}\label{parameta-orientation}
Let $C=(c_{ij}) \in M_n(\mathbb{Z})$ be a symmetrizable GCM and $D=\diag{c_1, \dots , c_n}$ be its symmetrizer.
\begin{enumerate2}
\item If $c_{ij} <0$, we set
		$g_{ij}\coloneqq |\gcdd{c_{ij}}{c_{ji}}|$ and $f_{ij}\coloneqq |c_{ij}| /g_{ij}{\,}$.
\item If $\Omega \subset \{1, 2, \dots , n\} \times \{1, 2, \dots , n\}$ satisfies the following two conditions, $\Omega$ is called an \textit{orientation} of $C$.
\begin{itemize2}
	\item If $\{(i, j), (j, i)\} \cap \Omega \neq\phi$ holds, then $c_{ij} < 0$;
	\item Every sequence $((i_1, i_2), (i_2, i_3), \dots , (i_t, i_{t+1})){\,\,}(t \geq 1)$ in $\{1, 2, \dots , n\} \times \{1, 2, \dots , n\}$ such that $(i_s, i_{s+1}) \in \Omega{\,}(1 \leq s \leq t)$ satisfies $i_1 \neq i_{t+1}$.
\end{itemize2}
    We define the opposite orientation ${\Omega}^* \coloneqq\{(j, i){\,}\mid(i, j) \in \Omega \}$ and $\overline{\Omega}\coloneqq\Omega \cup {\Omega}^*$.
\end{enumerate2}
\end{definition}
\begin{definition} \label{C-quiver}
	Under the setting of Definitions \ref{GCM}, \ref{parameta-orientation}, we define:
	\begin{enumerate2}
	\item The quiver $Q= Q{\,}(C, \Omega)= (Q_0, Q_1, s, t)$:
	\begin{align*}
		&Q_0= \{1, 2, \dots , n\},\\
		&Q_1= \{{\alpha}^{(g)}_{ij} \colon j \rightarrow i{\,}\mid(i, j) \in \Omega, 1 \leq g \leq g_{ij} \} \cup \{ {\varepsilon}_i \colon i \rightarrow i \mid 1 \leq i \leq n \}\\
		&s({\alpha}^{(g)}_{ij}) =j,{\,} t({\alpha}^{(g)}_{ij})=i,{\,} s({\varepsilon}_i) =t({\varepsilon}_i) =i\text{.}
	\end{align*}
	\item The \textit{double quiver} $\overline{Q} =(\overline{Q}_0, \overline{Q}_1, s, t)$ of $Q$ as follows:
	\begin{align*}
		&\overline{Q}_0= Q_0= \{1, 2, \dots , n\},\\
		&\overline{Q}_1= \{{\alpha}^{(g)}_{ij} \colon j \rightarrow i{\,}\mid(i, j) \in \overline{\Omega}, 1 \leq g \leq g_{ij} \} 
		\cup \{ {\varepsilon}_i \colon i \rightarrow i \mid 1 \leq i \leq n \},\\
		&s({\alpha}^{(g)}_{ij}) =j,{\,} t({\alpha}^{(g)}_{ij})=i,{\,} s({\varepsilon}_i) =t({\varepsilon}_i) =i\text{.}
	\end{align*}
\end{enumerate2}	
\end{definition}
Finally, we define $K$-algebras $H$ and $\Pi$ as quiver with relations.
\begin{definition} \label{def of H&Pi}
Under the setting of Definition \ref{C-quiver}, we define:
\begin{enumerate2}
\item[(1).] We define a $K$-algebra $H= H(C, D, \Omega)\coloneqq KQ/I$ by the quiver $Q$ with relations $I$ generated by (H1), (H2):
		\begin{enumerate2}
			\item[(H1)] $\eps{i}^{c_i}=0$\,($i\in Q_0$);
			\item[(H2)] For each $(i, j) \in \Omega$, we have $\eps{i}^{f_{ji}}\alpha^{(g)}_{ij}=\alpha^{(g)}_{ij}\eps{j}^{f_{ij}}$\, ($1 \leq g \leq g_{ij}$).
		\end{enumerate2}
	\item[(2).]	We define a $K$-algebra $\Pi= \fgpre{C}{D, \overline{\Omega}} \coloneqq K\overline{Q}/\overline{I}$ by the quiver $\overline{Q}$ with relations $\overline{I}$ generated by (P1)-(P3):
		\begin{enumerate2}
			\item[(P1)] $\eps{i}^{c_i}=0$ ($i\in Q_0$);
			\item[(P2)] For each $(i, j) \in \overline{\Omega}$, we have $\eps{i}^{f_{ji}}\alpha^{(g)}_{ij}=\alpha^{(g)}_{ij}\eps{j}^{f_{ij}}$ ($1 \leq g \leq g_{ij}$);
			\item[(P3)] For each $i \in Q_0$, we have
			$$\sum_{j \in \overline{\Omega}{\,}(-, i)} \sum_{g=1}^{g_{ij}} \sum_{f=0}^{{f_{ji}}-1} \sgn{i}{j}\eps{i}^f \alpha^{(g)}_{ij} \alpha^{(g)}_{ji} \eps{i}^{f_{ji}-1-f}=0, $$
			where we define
			\begin{align*}
	\overline{\Omega}{\,}(i, -)\coloneqq \{j\in Q_0\mid (i, j)\in \overline{\Omega}\}, &&\overline{\Omega}{\,}(-, j)\coloneqq\{i\in Q_0\mid (i, j)\in \overline{\Omega}\},
\end{align*}
and 
			$$\sgn{i}{j}\coloneqq 
			 \begin{cases}
				1 &(i, j)\in \Omega,\\
				-1 &(i, j)\in \Omega^*.
			\end{cases}$$
		\end{enumerate2}
		\end{enumerate2}
	We refer to the above $\Pi$ as the \textit{generalized preprojective algebra} associated with the pair $(C, D)$. In the definition of $\Pi$, let $\{ e_i\mid i\in Q_0\}$ be the complete set of primitive orthogonal idempotents corresponding to the vertex set $Q_0$. We note that $\Pi=\fgpre{C}{D}$ does not depend on a choice of orientations $\Omega$ up to isomorphism.
\end{definition}
\begin{example} \label{exampreproj}
	Let $C=\left(
	\begin{array}{rr}
		2 & -1 \\
	   -2 & 2 \\
	\end{array}
	\right)
	$, 
	$D=\diag{2d, d}\,(d\in \mathbb{Z}_{>0})$ and $\Omega=\{(1, 2)\}$. We have $c_1=2d, c_2=d, g_{12}=g_{21}=1, f_{12}=1, f_{21}=2$. Then, $\Pi=\fgpre{C}{D}$ is isomorphic to the $K$-algebra defined as the quiver
	$$
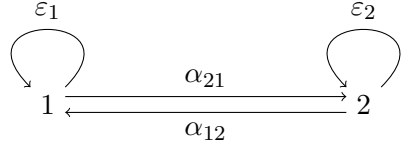
\begin{tikzpicture}[auto]
\node (a) at (0, 0) {$1$}; \node (b) at (4.2, 0) {$2$};
\draw[->, loop] (a) to node[swap] {$\eps{1}$} (a);
\draw[->, transform canvas={yshift=3pt}] (a) to node {$\alpha_{21}$} (b);
\draw[->, transform canvas={yshift=-3pt}] (b) to node {$\alpha_{12}$} (a);
\draw[->, loop] (b) to node[swap] {$\eps{2}$} (b);
\end{tikzpicture}$$
	with relations {(P1)} $\eps{1}^{2d}=0, \eps{2}^d=0$; {(P2)} $\eps{1}^2\alpha_{12}=\alpha_{12}\eps{2}, \eps{2}\alpha_{21}=\alpha_{21}\eps{1}^2$; {(P3)} $\alpha_{12}\alpha_{21}\eps{1}+\eps{1}\alpha_{12}\alpha_{21}=0, -\alpha_{21}\alpha_{12}=0$.
\end{example}
As a $K$-algebra, preprojective algebra of Dynkin type is characterized by the following Proposition.
\begin{proposition}[{Geiss-Leclerc-Schr\"oer \cite[Corollary 11.3, 12.7]{MR3660306}}]  \label{Piselfinj}
	Let $C$ be a connected GCM. Then, $\fgpre{C}{D}$ is a finite dimensional self-injective $K$-algebra if and only if $C$ is of Dynkin type.
\end{proposition}
\begin{definition}[Locally free modules] \label{def of H_i&lfrep}
	Under the setting of Definition \ref{def of H&Pi}, we define:
	\begin{enumerate2}
		\item $H_i \coloneqq e_i H e_i \cong  K[\eps{i}]/(\eps{i}^{c_i})$ for each $i \in Q_0$.
		\item A $\Pi$-module $M$ is \textit{locally free}, if $e_i M$ is a free $H_i$-module for each $i\in Q_0$.	
\end{enumerate2}
\end{definition}
We denote the full subcategory consisting of locally free modules (resp. locally free modules such that each free $H_i$-module $e_iM$ is of finite rank) by $\lfRep{\Pi}$ (resp. $\lfrep{\Pi}$).
\begin{definition}
	For each $i\in Q_0$, we say that $E_i\in\lfrep{\Pi}$ is a \textit{generalized simple module}, if we have the $H_i$-module isomorphisms,
	$$e_j E_i \cong\begin{cases}
		H_i &(j=i)\\
		0 &(j\neq i).
	\end{cases}$$
\end{definition}
By this definition, $E_i$ is a uniserial module that has only simple modules $S_i$ as the composition factors. We can also define locally free $\op{\Pi}$-modules and generalized simple $\op{\Pi}$-modules similarly. As $\Pi$ is a finite dimensional $K$-algebra in case $C$ is of Dynkin type, we have the standard $K$-duality $\D{(-)}\coloneqq\Hom{K}{-}{K}\colon\rep{\Pi}\longrightarrow\rep{\op{\Pi}}$. By definition of locally free modules, $M\in\lfrep{\Pi}$ if and only if $\D{M}\in\lfrep{\op{\Pi}}$. We know the following properties about locally free modules.
\begin{proposition}[{\cite[Lemma 3.8]{MR3660306}, Fu-Geng \cite[Lemma 2.6]{Fu2018}}] \label{rep-closed}
	$\lfRep{\Pi}$ is closed under kernel of epimorphisms, cokernel of monomorphisms, and extensions. 
\end{proposition}
\begin{proposition}[{\cite[Corollary 2.7]{Fu2018}}] \label{Pi_projdim}
	Let be $M\in\Rep{\Pi}$. If $\pdim{\Pi}{M}<\infty$, then $M\in \lfRep{\Pi}$.
\end{proposition}
\begin{proposition}[{\cite[Theorem 12.6]{MR3660306}}]\label{PiExtduality}
	For $M\in\lfRep{\Pi}$ and $N\in\lfrep{\Pi}$, we have the following functorial isomorphisms:
\begin{enumerate2}		
\item $\Ext{1}{\Pi}{M}{N} \cong \D{\Ext{1}{\Pi}{N}{M}}$.
\item If $C$ does not contain any components of Dynkin type, then we have $\Ext{2-i}{\Pi}{M}{N} \cong \D{\Ext{i}{\Pi}{N}{M}}$ for $i=0, 1, 2$.
\end{enumerate2}
\end{proposition}
We note that above three propositions hold for $\Rep{\op{\Pi}}$ and $\rep{\op{\Pi}}$.
\subsection{$\tau$-tilting theory}\label{ttilt}
In this subsection, we review the $\tau$-tilting theory due to Adachi-Iyama-Reiten \cite{MR3187626}. The basic references are \cite{MR3187626} and Demonet-Iyama-Jasso \cite{Demonet2015}. Let $\Lambda$ be a basic finite dimensional $K$-algebra. Let $\tau$ and $\tau^-$ be the Auslander-Reiten translations for $\rep{\Lambda}$, which give equivalences
$\tau \colon \pstrep{\Lambda} \rightarrow \istrep{\Lambda}$
    and
    $\tau^- \colon \istrep{\Lambda} \rightarrow \pstrep{\Lambda}$
between the projectively stable category $\pstrep{\Lambda}$ and the injectively stable category $\istrep{\Lambda}$. We denote the number of non-isomorphic indecomposable direct summand of $M$ by $\summand{M}$ and the two-sided ideal generated by an element $e\in\Lambda$ by $\langle e\rangle$.
\begin{definition} \label{ttilt-theory}
	Let $M\in\rep{\Lambda}$ and $P\in\proj{\Lambda}$. We define:
	\begin{enumerate2}
		\item $M$ is a \textit{$\tau$-rigid} $\Lambda$-module, if $\Hom{\Lambda}{M}{\tau M} =0$;
		\item $M$ is a \textit{$\tau$-tilting} $\Lambda$-module, if $M$ is $\tau$-rigid and $\summand{M} = \summand{\Lambda}$;
		\item $M$ is an \textit{almost complete $\tau$-tilting} $\Lambda$-module, if $M$ is $\tau$-rigid and $\summand{M} = \summand{\Lambda}-1$;
		\item $M$ is a \textit{support $\tau$-tilting} module, if there exists an idempotent $e \in \Lambda$ such that $M$ is $\tau$-tilting $(\Lambda / \langle e \rangle )$-module;
		\item $(M, P)$ is a \textit{$\tau$-rigid pair}, if $M$ is $\tau$-rigid and $\Hom{\Lambda}{P}{M} =0$;
		\item $(M, P)$ is a \textit{support $\tau$-tilting pair} (resp. \textit{almost complete $\tau$-tilting pair}), if $(M, P)$ is a $\tau$-rigid pair and $\summand{M} + \summand{P} = \summand{\Lambda}$ (resp. $\summand{M} + \summand{P} = \summand{\Lambda}-1$);
		\item $(M, P)$ is a \textit{direct summand of $(M', P')$}, if $(M, P)$ and $(M', P')$ are $\tau$-rigid pair and $M$ (resp. $P$) is a direct summand of $M'$  (resp. a direct summand of $P'$);
		\item $(M, P)$ is \textit{basic}, if $M$ and $P$ are basic. (i.e. each direct summand of $M\oplus P$ is multiplicity free).
	\end{enumerate2}
	\end{definition}
	We denote the full subcategory of indecomposable $\tau$-rigid $\Lambda$-modules by $\itrigid{\Lambda}$ and the full subcategory of basic support $\tau$-tilting $\Lambda$-modules by $\sttilt{\Lambda}$. We can think of $\tau$-rigid modules as a generalization of classical partial tilting modules in the sense of classical Bongartz's lemma by the following theorem:
	\begin{theorem}[{\cite[Theorem 2.10]{MR3187626}}] \label{Bongartz}	Any $\tau$-rigid $\Lambda$-module is a direct summand of some $\tau$-tilting $\Lambda$-module.
	\end{theorem}
	We denote the full subcategory of finite direct summands of finite direct sums of $M$ by $\add{M}$. A characterization of $\tau$-rigid pairs and support $\tau$-tilting pairs is given in the following theorem:
	\begin{theorem}[{\cite[Proposition 2.3]{MR3187626}}] \label{tau-pair}	 Let $M \in \rep{\Lambda}$, $P \in \proj{\Lambda}$ and $e \in \Lambda$ be an idempotent such that $\add{P} =\add{\Lambda e}$. 
	\begin{enumerate2}
		\item $(M, P)$ is a $\tau$-rigid pair, if and only if $M$ is a $\tau$-rigid $(\Lambda /\langle e \rangle)$-module;
		\item  $(M, P)$ is a support $\tau$-tilting pair if and only if $M$ is a $\tau$-tilting $(\Lambda /\langle e \rangle)$-module;
		\item $(M, P)$ is an almost complete support $\tau$-tilting pair, if and only if $M$ is an almost complete $\tau$-tilting $(\Lambda /\langle e \rangle)$ module;
		\item If $(M, P)$ and $(M, Q)$ are support $\tau$-tilting pairs in $\rep{\Lambda}$, then $\add{P} = \add{Q}$.
	\end{enumerate2}
\end{theorem}
By Theorem \ref{tau-pair}, we can identify basic support $\tau$-tilting modules with basic support $\tau$-tilting pairs. We define dual notions of $\tau$-rigid, $\tau$-tilting and support $\tau$-tilting modules.
\begin{definition}
Let $M \in \rep{\Lambda}$. We define:
    	\begin{enumerate2}
		\item $M$ is a \textit{$\tau^-$-rigid} $\Lambda$-module, if $\Hom{\Lambda}{\tau^- M}{M} =0$;
		\item $M$ is a \textit{$\tau^-$-tilting} $\Lambda$-module, if $M$ is $\tau^-$-rigid and $\summand{M} = \summand{\Lambda}$;
		\item $M$ is a \textit{support $\tau^-$-tilting} module, if there exists an idempotent $e \in \Lambda$ such that $M$ is $\tau^-$-tilting $(\Lambda / \langle e \rangle )$-module.
		\end{enumerate2}
\end{definition}
Note that $M$ is a $\tau^-$-rigid (resp. $\tau^-$-tilting, support $\tau^-$-tilting) $\Lambda$-module if and only if $\D M$ is a $\tau$-rigid (resp. $\tau$-tilting, support $\tau$-tilting) $\op{\Lambda}$-module by definition. We denote the set of basic support $\tau^-$-tilting modules by $\stmtilt{\Lambda}$.
\begin{definition}[{cf.~\cite[Proposition 1.1]{MR3187626}}]
\begin{enumerate2}
	\item A full subcategory $\mathcal{T}$ in $\rep{\Lambda}$ (resp. $\mathcal{F}$ in $\rep{\Lambda}$) is a \textit{torsion class} (resp. a \textit{torsion-free class}), if $\mathcal{T}$ (resp. $\mathcal{F}$) is closed under extensions and taking a factor module of objects (resp. taking a submodule of objects).
	\item A torsion class $\mathcal{T}$ in $\rep{\Lambda}$ (resp. a torsion-free classes $\mathcal{F}$ in $\rep{\Lambda}$) is \textit{functorially finite}, if there exists $M \in \rep{\Lambda}$ such that $\mathcal{T} =\Fac{M}$ (resp. $\mathcal{F} =\Sub{M}$), where $\Fac{M}$ (resp. $\Sub{M}$) is the full subcategory of factor modules (resp. submodules) of finite direct sums of $M$ in $\rep{\Lambda}$.
\end{enumerate2}
\end{definition}
We denote the set of torsion classes in $\rep{\Lambda}$  (resp. torsion-free classes) by $\tors{\Lambda}$ (resp. $\torf{\Lambda}$) and the set of functorially finite torsion classes in $\rep{\Lambda}$ (resp. torsion-free classes) by $\ftors{\Lambda}$ (resp. $\ftorf{\Lambda}$). In the $\tau$-tilting theory, one of the most important classes of algebras is $\tau$-tilting finite algebras:
\begin{def-pro}[cf.~\cite{Iyama2016}, \cite{Demonet2015}]
	An algebra $\Lambda$ is called \textit{$\tau$-tilting finite}, if $\Lambda$ satisfies one of the following equivalent conditions:
	\begin{itemize2}\label{ttilt-finite}
	    \item There are only finitely many isoclasses of basic $\tau$-tilting modules;
		\item $\sttilt{\Lambda}$ is a finite set;
		\item $\itrigid{\Lambda}$ is a finite set;
		\item $\ftors{\Lambda}$ (resp. $\ftorf{\Lambda}$) is a finite set;
		\item The poset $(\ftors{\Lambda}, \subseteq)$ (resp. $(\ftorf{\Lambda}, \subseteq)$) forms a complete lattice;
		\item $\ftors{\Lambda} = \tors{\Lambda}$;
		\item $\ftorf{\Lambda} = \torf{\Lambda}$.
	\end{itemize2}
\end{def-pro}
Let $\Lambda$ be a finite dimensional $K$-algebra again.
\begin{theorem}[{\cite[Theorem 2.7, 2.15]{MR3187626}}] \label{sttilt-tors} 
\begin{enumerate2}
\item We have the bijection between $\sttilt{\Lambda}$ and $\ftors{\Lambda}$:
	\begin{align*}
	\sttilt{\Lambda} &\longrightarrow \ftors{\Lambda}\\
	M &\longmapsto \Fac{M} \text{.}
	\end{align*}
\item We have the bijection between $\stmtilt{\Lambda}$ and $\ftorf{\Lambda}$:
	\begin{align*}
	\stmtilt{\Lambda} &\longrightarrow \ftorf{\Lambda}\\
	M &\longmapsto \Sub{M} \text{.}
	\end{align*}
\end{enumerate2}
\end{theorem}
We can define a partial order of $\sttilt{\Lambda}$ by Theorem \ref{sttilt-tors} as follows: 
\begin{definition}\label{poofsttilt}
	For $T, T' \in \sttilt{\Lambda}$, we define a partial order $\leq$ on $\sttilt{\Lambda}$ by $T \leq T' \Leftrightarrow \Fac{T} \subseteq \Fac{T'}$.\end{definition}
	Finally, the above partial order is understood in the terms of mutations:
	\begin{def-pro}[{\cite[Theorem 2.18]{MR3187626}}]\label{defofmutation}
		Any basic almost complete $\tau$-tilting pair $(U, Q)$ is a direct summand of precisely two different basic support $\tau$-tilting pairs $(T, P)$ and $(T', P')$. In addition, these $T, T'\in \sttilt{\Lambda}$ satisfy $T'<T$ or $T'>T$. In this setting, if $T'<T$ (resp. $T'>T$), we say that $(T', P')$ is a \textit{left} (resp. \textit{right}) \textit{mutation} of $(T, P)$. For this $T\in \sttilt{\Lambda}$ and the indecomposable summand $X$ of $T$ such that $T=X\oplus U$, we say that $T'$ is the \textit{left} (resp. \textit{right}) \textit{mutation} of $T$ at $X$, if $T'<T$ (resp. $T'>T$).
	\end{def-pro}
	\begin{proposition}[{\cite[Definition-Proposition 2.28]{MR3187626}}]\label{equivcondofmutation}
	Under the setting of Definition-Proposition \ref{defofmutation}, $T'$ is the left mutation of $T$ at $X$ if and only if $X\notin \Fac{U}$.
	\end{proposition}
	\begin{theorem}[{\cite[Theorem 2.33]{MR3187626}}]\label{mutation&order}
	Let $T, U \in\sttilt{\Lambda}$. The following conditions are equivalent:
	\begin{enumerate2}
		\item $U$ is a left mutation of $T$;
		\item $T$ is a right mutation of $U$;
		\item $T$ and $U$ satisfy $T>U$, and there does not exist $V\in\sttilt{\Lambda}$ such that $T>V>U$.
	\end{enumerate2}
\end{theorem}
We know that for a $\tau$-tilting finite algebra $\Lambda$, $\sttilt{\Lambda}$ forms a finite complete lattice by Definition-Proposition \ref{ttilt-finite} and Definition \ref{poofsttilt}. In particular, we find that $M, N \in \sttilt{\Lambda}$ are related by a mutation if and only if one is next to the other in the finite complete lattice of $\sttilt{\Lambda}$ by Theorem \ref{mutation&order}. Finally, we review a characterization of $\itrigid{\Lambda}$ in terms of the lattice.
\begin{theorem}[{\cite[Theorem 2.7 and its proof]{Iyama2016}}]\label{trigid-jirrtors}
	Let $\Lambda$ be a $\tau$-tilting finite algebra. Then, we have the following bijection:
	\begin{align*}
		\itrigid{\Lambda} &\longrightarrow \jirr{(\tors{\Lambda})}\\
		L &\longmapsto \Fac{L}.
	\end{align*}
	The inverse map is given by $\mathcal{T} \longmapsto N$, where $N$ is a unique indecomposable summand of $M\in\sttilt{\Lambda}$ giving $\mathcal{T}=\Fac{M}$ such that $\Fac{N} = \Fac{M}$. Now, $M\in\sttilt{\Lambda}$ has a unique indecomposable summand such that $\Fac{N} = \Fac{M}$ if and only if $M\in\sttilt{\Lambda}$ has a unique indecomposable summand $N$ such that $N \notin \Fac{M/N}$, equivalently $M$ has a unique left mutation in $\sttilt{\Lambda}$.
\end{theorem}
\subsection{Idempotent two-sided ideals of preprojective algebras and Weyl groups}\label{idem}
In this subsection, we review the work of Fu-Geng \cite{Fu2018} about a relationship between generalized preprojective algebras and Weyl groups of Kac-Moody Lie algebras. In \cite{Fu2018}, $K$ is assumed to be an algebraically closed field, but the discussion in \cite{Fu2018} works in the situation that $K$ is an arbitrary field. The basic materials about Coxeter groups appeared in this paper are found in \cite{MR2133266}. Let $\Pi =\fgpre{C}{D}$ be a generalized preprojective algebra associated with a symmetrizable GCM $C$ and its symmetrizer $D$. Let $W(C)$ be the Weyl group of the Kac-Moody Lie algebra associated with $C$. We define the \textit{idempotent ideal} $I_i$ to $i\in Q_0$ by $I_i\coloneqq \Pi (1-e_i)\Pi$. In particular, the inclusion $I_i\subseteq \Pi$ induces an exact sequence:
$$0\longrightarrow I_i \longrightarrow \Pi \longrightarrow E_i\longrightarrow 0.$$
\begin{theorem}[{\cite[Theorem 4.7]{Fu2018}}]\label{WandI_w}
	Let $C\in M_n(\mathbb{Z})$ be a symmetrizable GCM and $D$ be any symmetrizer $D$ of $C$. There is a bijection $\psi$ from $W(C)$ to the monoid $\langle I_1, I_2, \dots, I_n \rangle\coloneqq\{I_{i_1} I_{i_2} \cdots I_{i_k}\mid i_1, i_2, \dots, i_k \in Q_0, k\geq0\}$ given by $$\psi{\,}(w) = I_w = I_{i_1} I_{i_2} \cdots I_{i_k} \quad (\text{$w = s_{i_1} s_{i_2} \cdots s_{i_k}$ is a reduced expression of $w\in W$}).$$
	Here, $\psi$ does not depend on a choice of reduced expressions of $w$.
\end{theorem}
In \cite{Fu2018}, the results of Buan-Iyama-Reiten-Scott \cite{MR2521253} and Mizuno \cite{Mizuno2014} are generalized for our situation.
\begin{definition}
	\begin{enumerate2}
		\item We call a two-sided ideal $T$ of $\Pi$ \textit{tilting ideal} if $T$ is a left tilting $\Pi$-module and a right tilting $\Pi$-module. 
		\item We call a two-sided ideal $T$ of $\Pi$ \textit{support $\tau$-tilting ideal} if $T$ is a left support $\tau$-tilting $\Pi$-module and a right support $\tau$-tilting $\Pi$-module.
	\end{enumerate2}
\end{definition}
\begin{theorem}[{\cite[Lemma 3.2, 3.9, Theorem 3.12, 5.14, 5.17]{Fu2018}}] \label{bijection_W_tilt}
	 Let $C\in M_n(\mathbb{Z})$ be a symmetrizable GCM and $D$ be any symmetrizer $D$ of $C$.
	\begin{enumerate2}
		\item If $C$ has no components of Dynkin type, $T \in \langle I_1, I_2, \dots, I_n \rangle$ if and only if $T$ is a cofinite tilting ideal of $\Pi$, where tilting is in the sense of Miyasita \cite{MR852914} or Happel \cite{MR935124}. In particular, any object in $\langle I_1, I_2, \dots, I_n \rangle \subseteq \Rep{\Pi}$ has projective dimension at most $1$.
		\item If $C$ is of Dynkin type, $T \in \langle I_1, I_2, \dots, I_n \rangle$ if and only if $T$ is a basic support $\tau$-tilting ideal. In particular, $\psi\colon w\mapsto I_w$ in Theorem \ref{bijection_W_tilt} gives a bijection between $W$ and $\sttilt{\Pi}$.
	\end{enumerate2}
	\end{theorem}
	We note that we can obtain the similar Theorem in $\Rep{\op{\Pi}}$ as Theorem \ref{bijection_W_tilt}. In the case that $C$ is of Dynkin type, we have a relationship between the right weak Bruhat order $\leq_R$ on $W(C)$ and the mutation in $\sttilt{\Pi}$.	From now on until the end of this subsection, let $C\in M_n(\mathbb{Z})$ be a symmetrizable GCM of Dynkin type and $D$ be any symmetrizer $D$ of $C$.
\begin{theorem}[{\cite[Lemma 5.11, Proposition 5.13 and its proof]{Fu2018}}] \label{Tleftmutation}
	Let $T\in\langle I_1, \dots , I_n \rangle$. If $TI_i\neq T$, then $T$ has a left mutation $TI_i$ at $Te_i$ in $\sttilt{\Pi}$.
\end{theorem}
\begin{theorem}[{\cite[Theorem 5.16]{Fu2018}}]\label{mutationI}
	For $i \in Q_0$ and $w \in W$, $I_w, I_{ws_i}\in \sttilt{\Pi}$ are related by a right or left mutation. In particular, if $\length{ws_i}>\length{w}$,
	$$I_{w s_i}=\begin{cases}
		I_w(1-e_i)\, &(I_wI_i=0)\\
		I_wI_ie_i\oplus I_w(1-e_i)\, &(I_wI_i\neq 0).
	\end{cases}$$
\end{theorem}
From the above, we give the following generalization of a result of Mizuno \cite{Mizuno2014}:
\begin{theorem}\label{poofsttiltofpi}
	Let $w\in W$ and $i\in Q_0$. The following are equivalent:
	\begin{enumerate2}
		\item $\length{w}<\length{ws_i}=\length{w}+1$;
		\item $I_w I_i \neq I_w$;
		\item $I_w$ has a left mutation $I_{ws_i}$ at $I_w e_i$.
	\end{enumerate2}
	\begin{proof}
		(1)$\Rightarrow$(2) follows from Theorem \ref{WandI_w}. (2)$\Rightarrow$(3) follows from Theorem \ref{Tleftmutation}. We show (3)$\Rightarrow$(1). We assume that $\length{w}>\length{ws_i}$ and put $u=ws_i$. Since $\length{u}<\length{us_i}$ and (1)$\Rightarrow$(2), we find that $I_u$ has a left mutation $I_uI_i=I_{us_i}=I_w$ at $I_ue_i$ by Theorem \ref{Tleftmutation}. Then, $I_u$ is a left and right mutation of $I_w$ at $I_we_i$. This yields that $I_we_i\in \Fac{I_w(1-e_i)}$ and $I_we_i\notin \Fac{I_w(1-e_i)}$ by Proposition \ref{equivcondofmutation}. This is a contradiction. Thus, (3)$\Rightarrow$(1) holds.
	\end{proof}
	\end{theorem}
	\begin{theorem}\label{ppoofsttiltofpi}
		Let $w\in W$ and $i\in Q_0$. The following are equivalent:
	\begin{enumerate2}
		\item $\length{w}>\length{ws_i}=\length{w}-1$;
		\item $I_w I_i = I_w$;
		\item $I_{ws_i}$ has a left mutation $I_w$ at $I_{ws_i}e_i$.
	\end{enumerate2}
	\end{theorem}
In particular, we find that $u \leq_R v$ in $W$ if and only if $I_u \geq I_v$ in $\sttilt{\Pi}$. That is, $(W, \leq_R)$ can be identified with $\op{(\sttilt{\Pi}, \leq )}$ as a poset. We note that we can consider the poset structure of $\sttilt{\op{\Pi}}$ similarly, if we consider the left weak order instead of a right weak Bruhat order.
\section{Support $\tau$-tilting ideals and locally free modules}\label{we}

In this section, let $C$ be a connected symmetrizable GCM of Dynkin type, $D$ be any symmetrizer of $C$, and $\Pi$ be the generalized preprojective algebra associated with $C$ and $D$. Then, the Weyl group $W\coloneqq W(C)$ for $C$ is a finite Coxeter groups, and $\Pi$ is a $\tau$-tilting finite algebra.
\subsection{Two-sided ideals $I_w$ and locally freeness}\label{twoside}

We have the following lemma from a classification of GCMs of Euclidean types:
\begin{lemma}[{cf.~Carter~\cite[Appendix]{MR2188930}}]\label{affineCartan}
	For any connected GCM $C=(c_{ij}) \in M_n(\mathbb{Z})$ of Dynkin type, there is a connected GCM of Euclidean type $\widetilde{C}=(\tilde{c}_{ij}) \in M_{n+1}(\mathbb{Z})$ $\,(i, j\in\{0, 1,\dots ,n\})$ and its symmetrizer $\widetilde{D}=\diag{\tilde{c}_0, \tilde{c}_1, \dots , \tilde{c}_n}$ such that $\tilde{c}_{ij}=c_{ij}$ and $\tilde{c}_k=c_k$ $(i,j,k=1,2,\dots , n)$.\qed
\end{lemma}
By Lemma \ref{affineCartan}, we can construct a generalized preprojective algebra $\fgpre{\widetilde{C}}{\widetilde{D}}$ from $C$ and $D$.
\begin{lemma}\label{affpreproj}
	For the generalized preprojective algebra $\Pi=\fgpre{C}{D}$ and the vertex set $Q_{0}=\{ 1, 2, \dots, n \}$, we denote the vertex set of the generalized preprojective algebra $\affPi=\fgpre{\widetilde{C}}{\widetilde{D}}$ by $\widetilde{Q}_0=\{0, 1, 2,\dots, n\}$. Then, we have an $K$-algebra isomorphism $\affPi/\langle \tilde{e}_0 \rangle \cong \Pi$.
	\begin{proof}
		We denote the generators of $\affPi=\fgpre{\widetilde{C}}{\widetilde{D}}$ by the symbols with tildes as like $\tilde{\eps{i}}$, in order to distinguish from the generators of $\Pi$. Then, we have the surjective algebra homomorphism $\pi\colon \affPi \longrightarrow \Pi$ defined by $\tilde{\alpha}^{(g)}_{ij}\mapsto \alpha^{(g)}_{ij}\, (i\neq j)$, $\tilde{\eps{i}}\mapsto\eps{i}\,(i\neq 0)$, $\tilde{e}_0\mapsto0$. In particular, $\Ker{\pi}=\langle \tilde{e}_0 \rangle$. So, we obtain $\affPi/\langle \tilde{e}_0 \rangle \cong \Pi$.
	\end{proof}
\end{lemma}
\begin{theorem}\label{I_w-locallyfree}
	For any generalized Cartan matrix $C$ and symmetrizer $D$, the two sided ideals $I_w$ of $\fgpre{C}{D}$ are locally free. In particular, any object in $\sttilt{\Pi}$ is locally free.
	\begin{proof}
		Let $\affPi$ be the generalized preprojective algebra $\fgpre{\widetilde{C}}{\widetilde{D}}$, where $(\widetilde{C}, \widetilde{D})$ is borrowed from Lemma \ref{affineCartan}. We denote generators or subsets in $\affPi$ by symbols with tilde as in Lemma \ref{affpreproj}. Let $\pi\colon\affPi\longrightarrow\Pi$ be the surjective homomorphism in Lemma \ref{affpreproj}. We can regard a $\Pi$-module as a $\affPi$-module via $\pi$. In particular, we have $\widetilde{H}_i\cong H_i$ for each $i\neq0$. Thus, any locally free $\Pi$-module $M$ can be seen as a locally free $\affPi$-modules such that $\tilde{e}_0M=0$. Now, $\Pi$ is a projective $\Pi$-module, so that this is a locally free $\Pi$-module by Proposition \ref{Pi_projdim}. Then, $\Pi$ is a locally free $\affPi$-module by the above discussion. Similarly, $\affPi$ is a locally free $\affPi$-module. So, we conclude that $\langle \tilde{e}_0\rangle$ is a locally free $\affPi$-module by the short exact sequence $0\rightarrow \langle \tilde{e}_0\rangle\rightarrow \affPi \rightarrow \affPi/\langle\tilde{e}_0\rangle\rightarrow 0$ in $\Rep{\affPi}$ and Proposition \ref{rep-closed}. We put $\affI{i}\coloneqq\affPi(1-\tilde{e}_i)\affPi$. We define the two-sided ideals $\affI{w}=\affI{i_1}\dotsm \affI{i_k}$ for a reduced expression $w=s_{i_1}\dotsm s_{i_k}\in W(C)$, by regarding $W(C)\subset\widetilde{W}=W(\widetilde{C})$ in a natural way. (This definition does not effect on reduced expressions). On the other hand, we have $\langle\tilde{e}_0\rangle\subseteq\affI{w}$ by $\langle \tilde{e}_0 \rangle \subseteq \affI{i_j}\,(j=1, \dots, k)$. We have the $\affPi$-module isomorphism, 
		$$\widehat{I}_w \slash \langle \tilde{e}_0 \rangle =(\widehat{I}_{i_1} \slash \langle \tilde{e}_0 \rangle)\cdots(\widehat{I}_{i_k} \slash \langle \tilde{e}_0 \rangle)\cong I_{i_1}\cdots I_{i_k}= I_w.$$
		Now, we regard the above isomorphism as a $\Pi$-module isomorphism. Since $\pdim{\affPi}{\affI{w}}<\infty$ by Theorem \ref{bijection_W_tilt} , $\affI{w}$ is a locally free $\affPi$-module. Then, $\affI{w}/\langle\tilde{e}_0\rangle$ is a locally free $\affPi$-module by Proposition \ref{rep-closed}. So, we find that $e_i I_w$ is a free $H_i$-module for each $i\in \{1, \dots, n\}$ via $\affPi/\langle \tilde{e}_0 \rangle \cong \Pi$ and $\widetilde{H}_i\cong H_i\,(i\neq0)$. This implies $I_w\in \lfrep{\Pi}$ for any $w\in W$.
	\end{proof}
\end{theorem}
As the corollary of Theorem \ref{I_w-locallyfree}, we have the following:
\begin{corollary} \label{locallyfreeitrigid}
	Any indecomposable $\tau$-rigid module in $\rep{\Pi}$ is locally free.
	\begin{proof}
		Since $\Pi$ is self-injective by Proposition \ref{Piselfinj}, all $\Pi e_i$ are indecomposable injective $\Pi$-modules. So, $\soc{\Pi e_i}$ are simple modules. In particular, the submodules $I_w e_i$ of $\Pi e_i$ are indecomposable. Since $\tau$-rigid modules appear as a summand of some $\tau$-tilting module by Theorem \ref{Bongartz}, any indecomposable $\tau$-rigid module is isomorphic to one of the modules $I_w e_i$ by Theorem \ref{bijection_W_tilt}. Since indecomposable summands of any locally free module is locally free, any indecomposable $\tau$-rigid $\Pi$-module is locally free.
	\end{proof}
\end{corollary}
We can classify $\tau$-rigid $\Pi$-modules in terms of the lattice of right weak Bruhat order on $W$.
\begin{theorem}\label{itrigidclass}
	The following map 
	\begin{align*}
		\mirr{W} &\longrightarrow \itrigid{\Pi}\\
		w &\longmapsto I_w e_k
		\end{align*}
		is bijective, where $k$ is the unique index such that $\length{ws_k} = \length{w} + 1$. 
		
		The proof of this Theorem is similar as the proof in \cite{Iyama2016}.
		\begin{proof}
			Since $\Pi$ is a $\tau$-tilting finite algebra, we have the bijective correspondence between $\sttilt{\Pi}$ and $\tors{\Pi}$
			\begin{align*}
		    \sttilt{\Pi} &\longrightarrow \tors{\Pi}\\
		    I_w &\longmapsto \Fac{I_w}.
		    \end{align*}
		    By definition of the order in $\sttilt{\Pi}$,  the posets $(W, \leq_R)$ and $\op{(\tors{\Pi}, \subseteq)}$ are isomorphic. That is, we have the following lattice isomorphism
		    \begin{align*}
		    	(W, \leq_R) &\longrightarrow \op{(\tors{\Pi}, \subseteq)}\\
		    	w &\longmapsto \Fac{I_w}.
		    \end{align*}
		    Thus, we have the following bijection between meet-irreducible elements of $W$ and join-irreducible elements of $\tors{\Pi}$
		    \begin{align*}
		    \mirr{W} &\longrightarrow \jirr{(\tors{\Pi})}\\
		    w &\longmapsto \Fac{I_w}.
		    \end{align*}
		    By Theorem \ref{trigid-jirrtors}, we have
		    \begin{align*}
		    \jirr{(\tors{\Pi})} &\longrightarrow \itrigid{\Pi}\\
		    \Fac{I_w} &\longmapsto I_w e_k,
		    \end{align*}
		    where $k$ is the unique index satisfying $\Fac{(I_w e_k)} = \Fac{I_w}$.
		    We have only to show that this $k$ is the unique index satisfying $\length{ws_k} > \length{w}$. Now, $I_{ws_k}$ and $I_w$ are objects  in a relation of mutation each other in $\sttilt{\Pi}$ by Theorem \ref{mutationI}.
		    In particular, if we give the decomposition $I_{ws_k}=\bigoplus_{p\in Q_0} (I_{ws_k}e_p)$ and $I_{w}=\bigoplus_{q\in Q_0} (I_{w}e_q)$, we have $\bigoplus_{p\neq k}I_{ws_k}e_p=\bigoplus_{q\neq k}I_{w}e_q$, and denote them by $U$. That is, $I_w = I_w e_k \oplus U$. Now, the condition that $k$ is the unique index such that $I_w e_k \notin \Fac(U)$ is equivalent to that $k$ is the unique index such that $\Fac{I_w s_k}=\Fac{I_w}$ by Theorem \ref{trigid-jirrtors}. So, since $k$ is the unique index such that $I_w s_k$ is a left mutation of $I_w$ at $I_w e_k$ by Proposition \ref{equivcondofmutation}, $k$ is the unique index such that $\length{ws_k}=\length{w}+1$ by Theorem \ref{poofsttiltofpi}. 
		\end{proof}
\end{theorem}
\begin{example}
In the type of $B_2$, the Hasse quivers of $(W, \leq_R)$ and $\op{(\sttilt{\Pi}, \leq)}$ are the following:
$$
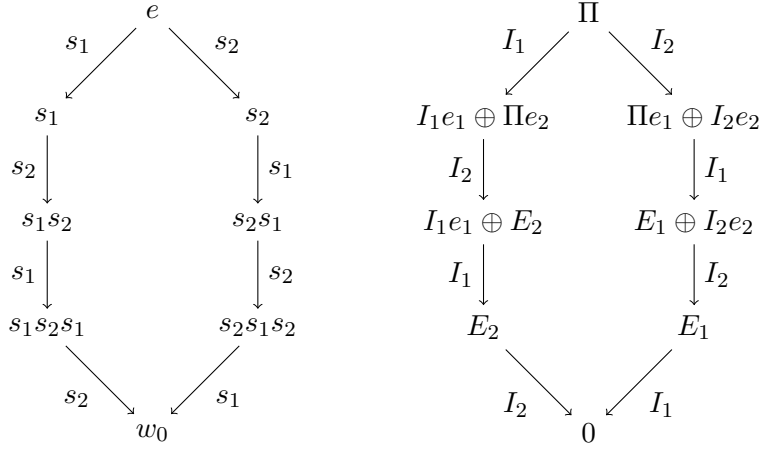
\begin{tikzpicture}[auto]
\node (a) at (1.4, 5.6) {$e$}; \node (b) at (0, 4.2) {$s_1$};
\node (c) at (2.8, 4.2) {$s_2$}; \node (d) at (0, 2.8) {$s_1s_2$};
\node (e) at (2.8, 2.8) {$s_2s_1$}; 
\node (f) at (0, 1.4) {$s_1s_2s_1$}; \node (g) at (2.8, 1.4) {$s_2s_1s_2$};
\node (h) at (1.4, 0) {$w_0$};
\node (1) at (7.2, 5.6) {$\Pi$}; \node (2) at (5.8, 4.2) {$I_1e_1\oplus\Pi e_2$};
\node (3) at (8.6, 4.2) {$\Pi e_1\oplus I_2e_2$}; \node (4) at (5.8, 2.8) {$I_1e_1\oplus E_2$};
\node (5) at (8.6, 2.8) {$E_1\oplus I_2e_2$}; 
\node (6) at (5.8, 1.4) {$E_2$}; \node (7) at (8.6, 1.4) {$E_1$};
\node (8) at (7.2, 0) {$0$}; 
\draw[->] (a) to node[swap] {$s_1$} (b);
\draw[->] (a) to node {$s_2$} (c);
\draw[->] (b) to node[swap] {$s_2$} (d);
\draw[->] (c) to node {$s_1$} (e);
\draw[->] (d) to node[swap] {$s_1$} (f);
\draw[->] (e) to node {$s_2$} (g);
\draw[->] (f) to node[swap] {$s_2$} (h);
\draw[->] (g) to node {$s_1$} (h);
\draw[->] (1) to node[swap] {$I_1$} (2);
\draw[->] (1) to node {$I_2$} (3);
\draw[->] (2) to node[swap] {$I_2$} (4);
\draw[->] (3) to node {$I_1$} (5);
\draw[->] (4) to node[swap] {$I_1$} (6);
\draw[->] (5) to node {$I_2$} (7);
\draw[->] (6) to node[swap] {$I_2$} (8);
\draw[->] (7) to node {$I_1$} (8);
\end{tikzpicture}.$$
In above picture, we have $\mirr{W}=\{$$s_1$, $s_2$, $s_1s_2$, $s_2s_1$, $s_1s_2s_1$, $s_2s_1s_2$$\}$ and $\itrigid{\Pi}=\{$$I_1e_2$, $I_2e_1$, $I_1I_2e_1$, $I_2I_1e_2$, $I_1I_2I_1e_2$, $I_2I_1I_2e_1$$\}=\{$$\Pi e_2$, $\Pi e_1$, $I_1e_1$, $I_2e_2$, $E_2$, $E_1$$\}$. 
\end{example}
\subsection{Torsion classes and torsion-free classes of $\rep{\Pi}$}\label{tor}
The following Lemma is shown in \cite{Fu2018} by using the minimal projective presentations of $\tau$-rigid $\Pi$-modules:
\begin{lemma}[{\cite[Lemma 5.9, 5.10]{Fu2018}}] \label{E_i tens or tor =0}
	For any $i\in Q_0$, we have:
	\begin{enumerate2}
\item For a support $\tau$-tilting $\Pi$-module $T$ and  a generalized simple module $E'_i$ in $\lfrep{\op{\Pi}}$, we have:
\begin{itemize2}
\item Either $E'_i \tenso{\Pi} T=0$ or $\Tor{1}{\Pi}{E'_i}{T}=0$;
\item $E'_i\tenso{\Pi} T=0$ if and only if $I_iT=T$;
\end{itemize2}
\item[\rm{(2).}] For a support $\tau$-tilting $\op{\Pi}$-module $T$ and a generalized simple module $E_i$ in $\lfrep{\Pi}$, we have:
\begin{itemize2}
\item Either $T \tenso{\Pi} E_i=0$ or $\Tor{1}{\Pi}{T}{E_i}=0$;
\item $T\tenso{\Pi} E_i=0$ if and only if $TI_i=T$.\qed
\end{itemize2}
\end{enumerate2}
\end{lemma}
We have:
\begin{proposition} \label{Extvan}
	For any $i\in Q_0$ and $w\in W$, we have the following statement:
	\begin{enumerate2}
		\item Let $I_w\in \sttilt{\Pi}$. If $\length{s_i w} > \length{w}$, then $\Ext{1}{\Pi}{I_w}{E_i}=0=\Ext{1}{\Pi}{E_i}{I_w}$.
		\item Let $I_w\in\sttilt{\op{\Pi}}$. If $\length{w s_i} > \length{w}$, then $\Ext{1}{\op{\Pi}}{I_w}{E'_i}=0=\Ext{1}{\op{\Pi}}{E'_i}{I_w}$.
		\end{enumerate2}
		\begin{proof}
			Since the proofs of the both cases are similar, we only prove (1). We have a short exact sequence:
			$$0 \rightarrow I_i \xrightarrow{\iota} \Pi \rightarrow E'_i \rightarrow 0$$
			in $\lfrep{\op{\Pi}}$, where $\iota\colon I_i\longrightarrow \Pi$ is the inclusion. Apply the functor $(-) \tenso{\Pi} I_w$ to this exact sequence to obtain the exact sequence:
			$$\Tor{1}{\Pi}{E'_i}{I_w} \longrightarrow I_i \tenso{\Pi} I_w \xrightarrow{\iota \tenso{\Pi} I_w} \Pi \tenso{\Pi} I_w \longrightarrow E'_i \tenso{\Pi} I_w \longrightarrow 0.$$
			If $\length{s_i w} > \length{w}$, $I_i I_w\subsetneq I_w$ by Theorem \ref{WandI_w}. So, $\Img{(\iota \tenso{\Pi} I_w)} =I_i I_w\subsetneq I_w$ via the natural isomorphism $\Pi \tenso{\Pi} I_w \cong I_w$, which implies $E'_i \tenso{\Pi} I_w \neq 0$. We obtain the following linear isomorphism,
			\begin{align*}
			\D{\Ext{1}{\Pi}{E_i}{I_w}}
			&\cong \Ext{1}{\Pi}{I_w}{E_i}&\text{($I_w \in \lfrep{\Pi}$ and  Proposition \ref{PiExtduality})}\\
			&\cong \Ext{1}{\Pi}{I_w}{\D{E'_i}}
			&(\D{E'_i}\cong E_i \in \lfrep{\Pi})\\
			&\cong \D{\Tor{1}{\Pi}{E'_i}{I_w}} &\text{(cf.~\cite[Appendix. Proposition 4.11]{MR2197389})}\\
			&=0&\text{(Lemma \ref{E_i tens or tor =0})}.
		\end{align*}
		This proves the assertion.
		\end{proof}
\end{proposition}
From the proof of Proposition \ref{Extvan}, we also deduce the following Proposition.
\begin{proposition}\label{IE-Tor-Tens}
	For any $i\in Q_0$ and $ w\in W$, we have:	\begin{enumerate2}
		\item Let $I_w\in\sttilt{\op{\Pi}}$. If $\length{ws_i} > \length{w}$, then $\Tor{1}{\Pi}{I_w}{E_i} = 0$;
		\item Let $I_w\in\sttilt{\Pi}$. If $\length{ s_i w} > \length{w}$, then $\Tor{1}{\Pi}{E'_i}{I_w} = 0$.\qed
	\end{enumerate2}
\end{proposition}
\begin{lemma} \label{E_i dsum}
	Let $i\in Q_0$.
	\begin{enumerate2}
		\item If $M, N \in \lfrep{\Pi}$ such that $M \subseteq N$, then $I_i M = I_i N$ if and only if there exists some non-negative integer $n$ such that $N/M \cong {E_i}^{\oplus n}$.
		\item If $M, N \in \lfrep{\op{\Pi}}$ and $M \subseteq N$, then $MI_i = NI_i$ if and only if  there exists some non-negative integer $n$ such that $N/M \cong {E'_i}^{\oplus n}$.
	\end{enumerate2}
	\begin{proof}
		The proof of both assertions are similar. So, we only prove (1). If $I_i M=I_i N$, then $I_i (N/M)=0$. Now, since $M, N \in \lfrep{\Pi}$, we have $N/M\in \lfrep{\Pi}$ by Proposition \ref{rep-closed}.  Now, since $I_i=\Pi(1-e_i)\Pi$, we have $(1-e_i)(N/M)=0$ and so $N/M=e_i (N/M)$. Thus, $N/M$ is isomorphic to a direct sum of some copies of the generalized simple module $E_i$ by definition of locally free modules.
		
		Conversely, if $N/M\cong E_i^{\oplus n}$, then we have $I_i \tenso{\Pi} (N/M) \cong I_i \tenso{\Pi} ({E_i}^{\oplus n})\cong (I_i \tenso{\Pi} {E_i})^{\oplus n}= 0$, because of Lemma \ref{E_i tens or tor =0} (2)-(ii) and the fact that $I_iI_i=I_i$. On the other hand, we have a short exact sequence:
		$$0\longrightarrow M\longrightarrow N\longrightarrow N/M\longrightarrow 0,$$
		and apply the functor $I_i\tenso{\Pi}(-)$ to obtain the exact sequence:
		$$I_i \tenso{\Pi} M \longrightarrow I_i \tenso{\Pi} N \longrightarrow 0.$$
		So, the inclusion $I_i M \hookrightarrow I_i N$ is surjective. That is, $I_i M = I_i N$ in $M\subseteq N$.
	\end{proof}
\end{lemma}
Finally, we obtain the following dualities. The proof is along the line of \cite{Oppermann2015}, but we need a little device for it.
\begin{theorem} \label{dualityofI}
 Let $w_0 \in W$ be the longest element in $W$. Then, we have the following isomorphisms.
 \begin{enumerate2}
		\item $\D{I_w} \cong \Pi / I_{w_0 w^{-1}}$ in $\lfrep{\op{\Pi}}$.
		\item $\D{I_w} \cong \Pi / I_{w^{-1} w_0}$ in $\lfrep{\Pi}$.
	\end{enumerate2}
	\begin{proof}
	We can prove (2) similarly as (1). So, we only prove (1). Since $\Pi$ is a self-injective algebra, we have $\Pi \cong \D{\Pi}$ in $\lfrep{\op{\Pi}}$. Applying $\D{(-)}$ to the inclusion $I_w\longrightarrow\Pi$ in $\rep{\Pi}$, we obtain an epimorphism $\D{\Pi}\longrightarrow \D{I_w}$ in $\lfrep{\op{\Pi}}$. Thus, we obtain an epimorphism $\psi\colon \Pi \longrightarrow \D{I_w}$ in $\lfrep{\op{\Pi}}$. We put $\tildeI{w}\coloneqq \Ker{\psi}$. Since $\tildeI{w}$ is a kernel of an epimorphism in $\lfrep{\op{\Pi}}$, we have $\tildeI{w}\in \lfrep{\op{\Pi}}$ by the right version of Proposition \ref{rep-closed}. It remains to show that $\tildeI{w}=I_{w_0 w^{-1}}$.
	
	Since $w_0\in W$ is the longest element, we obtain $I_{w_0w^{-1}}I_w=I_{w_0}=0$, and so $I_{w_0w^{-1}}$ is contained in the right annihilator of $\D{I_w}$. That is, $\tildeI{w}\supseteq I_{w_0w^{-1}}$.
	
	We show the converse inclusion $\tildeI{w}\subseteq I_{w_0w^{-1}}$ by induction on $\length{w}$. When $w$ is an identity, $\widetilde{I}_w =0=I_{w_0 w^{-1}}$ because $I_w=\Pi$ and $I_{w_0}=0$. So, we assume that $\length{s_i w} > \length{w}$. Now, $I_w /I_{s_i w}$ is the cokernel of the inclusion $I_{s_i w} = I_i I_w \subsetneq I_w$ in $\lfrep{\Pi}$. Since $I_iI_i=I_i$ and so $I_iI_{s_iw}=I_iI_w$, there exists some non-negative integer $n$ such that $I_w /I_{s_i w}\cong E_i^{\oplus n}$ by Lemma \ref{E_i dsum} (1). Thus, we have an exact sequence:
	$$0\longrightarrow I_{s_iw}\longrightarrow I_w\longrightarrow E_i^{\oplus n}\longrightarrow 0$$
	in $\lfrep{\Pi}$. Applying $\D{(-)}$ yields an exact sequence:
	$$0 \longrightarrow {E'_i}^{\oplus n} \longrightarrow \D{I_w} \longrightarrow \D{I_{s_i w}} \longrightarrow 0$$
	in $\lfrep{\op{\Pi}}$. On the other hands, the inclusion $\tildeI{w} \hookrightarrow \tildeI{s_iw}$ gives the short exact sequence:
	$$0 \longrightarrow \tildeI{w} \longrightarrow \tildeI{s_i w} \longrightarrow \tildeI{s_i w}/\tildeI{w} \longrightarrow 0.$$
	in $\lfrep{\op{\Pi}}$. So, we obtain the following commutative diagram with exact rows in $\lfrep{\op{\Pi}}$:
	$$   	\begin{tikzpicture}[auto]
\node (a) at (0, 0) {$0$}; \node (x) at (1.4, 0) {$\tildeI{s_i w}$};
\node (b) at (0, 1.4) {$0$};   \node (y) at (1.4, 1.4) {$\tildeI{w}$};
\node (z) at (1.4, -1.4) {$\tildeI{s_i w}/\tildeI{w}$}; 
\node (c) at (2.8, 0) {$\Pi$}; \node (d) at (2.8, 1.4) {$\Pi$};
\node (e) at (4.2, 0) {$\D{I_{s_i w}}$}; \node (f) at (4.2, 1.4) {$\D{I_w}$}; \node (g) at (4.2, 2.8) {${E'_i}^{\oplus n}$};
\node (h) at (5.6, 0) {$0$}; \node (i) at (5.6, 1.4) {$0$};
\draw[->] (a) to (x);
\draw[right hook->] (y) to (x);
\draw[->] (y) to (d);
\draw[->] (b) to (y);
\draw[->] (x) to (c);
\draw[->>] (x) to (z);
\draw[->] (d) to node {$\mathrm{id}_{\Pi}$} (c);
\draw[->] (d) to (f);
\draw[->] (c) to (e);
\draw[->>] (f) to (e);
\draw[->] (f) to (i);
\draw[->] (e) to (h);
\draw[right hook->] (g) to (f);
\end{tikzpicture}.$$
We have $\tildeI{s_i w}/\tildeI{w} \cong {E'_i}^{\oplus n}$ by the snake lemma. It follows that $\tildeI{s_i w} I_i = \tildeI{w} I_i$ by Lemma \ref{E_i dsum}, and hence $\tildeI{s_i w} I_i = \tildeI{w} I_i\subseteq \tildeI{w} = I_{w_0 w^{-1}}$ by induction hypothesis. Now, the inequality $\length{w_0 w^{-1} s_i} = \length{w_0} - \length{w^{-1} s_i}
		= \length{w_0} - \length{s_i w}
		< \length{w_0} - \length{w}
		= \length{w_0 w^{-1}}$ gives $I_{w_0 w^{-1}} = I_{w_0 w^{-1} s_i} I_i$ for $I_{w_0w^{-1}} \in \sttilt{\op{\Pi}}$. So, we have $\tildeI{s_i w} I_i\subseteq I_{w_0 w^{-1} s_i} I_i$. On the other hand, we have $I_{w_0 w^{-1} s_i} \subseteq \tildeI{s_i w}$ in the first part of this proof, and so $I_{w_0 w^{-1} s_i} I_i \subseteq \tildeI{s_i w} I_i$. Thus, we obtain $\tildeI{s_i w} I_i = I_{w_0 w^{-1} s_i} I_i$.
		
		By the right version of Theorem \ref{I_w-locallyfree} and Lemma \ref{E_i dsum}, there exists non-positive integer $m$ such that $\tildeI{s_i w}/ I_{w_0 w^{-1} s_i} \cong {E'_i}^{\oplus m}$. So, we only have to show that $m=0$. Now, since $\length{w_0 w^{-1} s_is_i}=\length{w_0 w^{-1}}>\length{w_0w^{-1} s_i}$, we have $\Ext{1}{\op{\Pi}}{E'_i}{I_{w_0 w^{-1} s_i}} =0$ by Proposition \ref{Extvan}. Thus, we obtain the following splitting exact sequence
		$$0 \longrightarrow I_{w_0 w^{-1} s_i} \longrightarrow \tildeI{s_i w} \longrightarrow {E'_i}^{\oplus m} \longrightarrow 0.$$
		So, there exists a right ideal $L$ of $\Pi$ isomorphic to ${E'_i}^{\oplus m}$ such that $\tildeI{s_i w} = L \oplus I_{w_0 w^{-1} s_i}$. Now, we assume that $L\neq 0$. Since $\Pi$ is a basic self-injective algebra, $\Pi$ decomposes into non-isomorphic indecomposable injective modules. That is, $\Pi =\bigoplus_{j\in Q_0}e_j\Pi$ is the direct sum of injective envelopes of simple $\op{\Pi}$-modules $S'_k\,(k\in Q_0)$. In particular, $\soc{\Pi}\cong\bigoplus_{k\in Q_0}S'_k$ and $S'_k\,(k\in Q_0)$ are non-isomorphic to each other. Since $E'_i$ is a uniserial module that has only $S'_i$ as the composition factors, we have $\soc{E'_i}\cong S'_i$. Thus, we obtain  $\soc{L}\cong {S'_i}^{\oplus m}$. Since $\soc{\Pi}\cong\bigoplus_{k\in Q_0}S'_k$ and $L$ is a right ideal of $\Pi$, we find $m=1$ and so $L \cong E'_i$. Now, we have $L = I_{w_0 s_i}$ by Lemma \ref{E_i dsum} and the general fact that for a basic finite dimensional self-injective $K$-algebra $\Lambda$ and a two-sided ideal $I$ of $\Lambda$, any right ideal $L$ of $\Lambda$ isomorphic to $I$ in $\rep{\op{\Lambda}}$ satisfies $L=I$ as sets (cf. \cite[proof of Lemma 2.20]{Mizuno2014}). In particular, we find that $\tildeI{s_i w} =  I_{w_0 s_i} \oplus I_{w_0 w^{-1} s_i}$ and so $I_{w_0 s_i} \cap I_{w_0 w^{-1} s_i} =0$. Now, since $w_0 s_i = w_0 w w_0 w_0 w^{-1} s_i$, we obtain the equality $\length{w_0 w w_0} + \length{w_0 w^{-1} s_i} =\length{w}+ \length{w_0} - \length{w^{-1}s_i}   =\length{w}+ \length{w_0} - \length{s_iw}=\length{w}+\length{w_0} - (1+\length{w})=\length{w_0}-1=\length{w_0 s_i}$. This equality gives $I_{w_0 s_i} = I_{w_0 w w_0} I_{w_0 w^{-1} s_i}$. However, this means that $I_{w_0 s_i} \subset I_{w_0 w^{-1} s_i}$, and so we obtain a contradiction. Thus, $L = 0$ and $\tildeI{s_i w} = I_{w_0 w^{-1} s_i}$. So, we conclude that $\tildeI{w} \subseteq I_{w_0 w^{-1}}$ for any $w \in W$ by induction on $\length{w}$ as required. 
		\end{proof}
\end{theorem}
As a corollary of Theorem \ref{dualityofI}, we obtain a classification theorem of torsion-free classes in $\rep{\Pi}$. Note that the former statement follows as a combination of \cite[Theorem 5.17]{Fu2018}, \cite[Theorem 2.7]{MR3187626} and \cite[Theorem 3.8]{Demonet2015}.
\begin{corollary} \label{Pitorstorf}
	We have the following two bijections:
	\begin{align*}
    W &\longrightarrow \tors{\Pi}
    &W &\longrightarrow \torf{\Pi}\\
    w &\longmapsto \Fac{I_w},
	&w &\longmapsto \Sub{(\Pi / I_w)}.
    \end{align*}
    \begin{proof}
    	The classification of torsion classes is a direct result of Theorem \ref{sttilt-tors} and Definition-Proposition \ref{ttilt-finite}.
    	Since $\sttilt{\op{\Pi}}$ corresponds to $\stmtilt{\Pi}$ bijectively via $I_w \mapsto \D{I_w}$, we have a bijection $W \rightarrow \stmtilt{\Pi}$ defined by $w \mapsto \Pi/I_w$ by Theorem \ref{dualityofI}. Thus, we obtain the second bijection by Theorem \ref{sttilt-tors} and Definition-Proposition \ref{ttilt-finite}.
    	    \end{proof}
\end{corollary} 
\begin{corollary}
	For $w\in W$, we have $\ann{I_w}=I_{w_0w^{-1}}$. In particular, the tilting objects in $\sttilt{\Pi}$ are nothing but $\Pi$ up to isomorphism.
	\begin{proof}
		Since the left annihilator of $I_w$ coincides with the right annihilator of $\D{I_w}$, the first assertion follows from Theorem \ref{dualityofI}.
		Since $w \mapsto I_w$ is a bijection from $W$ to $\sttilt{\Pi}$ by Theorem \ref{bijection_W_tilt}, the object $I_{w} \in \sttilt{\Pi}$ satisfying $\ann{I_w} = I_{w_0 w^{-1}} =0$ are nothing but $I_{e} = \Pi$ for the identity $e \in W$. Then, the second assertion follows from the general fact \cite[Proposition 2.2]{MR3187626} of $\tau$-tilting theory that faithful support $\tau$-tilting modules are precisely tilting modules in the sense of Brenner-Butler \cite{MR607151}.
	\end{proof}
\end{corollary}
\subsection*{Acknowledgement}
The author thanks his supervisor Syu Kato for hopeful encouragement and pointing out many typographical errors of a draft of this paper. The author also thanks Bernard Leclerc for telling him an easy proof of Theorem \ref{locallyfreeitrigid} and giving him many interesting lectures about this topic during his stay in Kyoto. Finally, the author thanks the referee for carefully reading the manuscript and for giving constructive comments.

\bibliographystyle{amsplain}
\bibliography{genprepro.bib}
\end{document}